\documentclass[12pt]{amsart}

\usepackage{amsmath,amsfonts,amssymb,amscd}
\begin{document}
\renewcommand{\thefootnote}{\fnsymbol{footnote}}
\pagestyle{plain}

\title{K-stability of constant scalar
curvature polarization}
\author{Toshiki Mabuchi${}^*$}
\maketitle
\abstract
In this paper, we shall show that a polarized algebraic manifold 
is K-stable if the polarization class
admits a K\"ahler metric of constant scalar curvature.
This generalizes the results of Chen-Tian \cite{CT},
Donaldson \cite{D3} and Stoppa \cite{S1}.
\endabstract
\section{Introduction}

Yau's conjecture${}^{\dagger}$ (\cite{Y1}, pp.49--50)
suggests a  strong correlation 
between stability of polarized algebraic manifolds and 
the existence of extremal metrics in the polarization class 
(cf. \cite{T1}, \cite{D0}, \cite{M1}).
Especially, for K\"ahler metrics of constant scalar curvature, 
the following is still open as an interesting conjecture:

\footnotetext{ ${}^{\dagger}$Its vector bundle counterpart known also as the Hitchin-Kobayashi correspondence 
is affirmative by the works of Kobayashi, L\"ubke, Donaldson, Uhlenbeck and Yau. }

 \medskip\noindent
{\bf Conjecture {\rm (Tian \cite{T1}, Donaldson \cite{D1})}.}
{\em A polarized algebraic manifold $(M,L)$ is K-stable if
and only if the polarization class $c_1(L)_{\Bbb R}$ admits a K\"ahler metric  
of constant scalar curvature.} 

\medskip
For ``if'' part of this conjecture,
Chen-Tian \cite{CT} and Donaldson \cite{D3}
showed that a polarized algebraic manifold $(M,L)$ 
is K-semistable when the class $c_1(L)_{\Bbb R}$ admits 
a K\"ahler metric of constant scalar curvature. 
Very recently, by an effective use of moduli spaces, Stoppa \cite{S1} proved a stronger result 
showing that a polarized algebraic manifold $(M,L)$ with a K\"ahler metric of constant scalar curvature  in $c_1(L)_{\Bbb R}$ is K-stable
if the group $\operatorname{Aut}(M,L)$ of holomorphic automorphisms 
of $(M,L)$ is discrete.  
The purpose of this paper is to extend Stoppa's result to the following general case
without assuming such discreteness:

\medskip \noindent
{\bf Main Theorem.} 
{\em A polarized algebraic manifold $(M,L)$ is K-stable if
the class $c_1(L)_{\Bbb R}$ admits a K\"ahler metric  
of constant scalar curvature.} 

\medskip
It should be emphasized that one of the main ingredients of this paper is the 
energy-theoretic approach to K-stability as in Tian \cite{T1} (see also \cite{CT}),
where in our actual proof of Main Theorem,  the logarithm of the Chow norm
is used in place of the K-energy. To see this clearly, 
for a connected Fano manifold $X$, we consider a {\it special degeneration} 
$$
\operatorname{pr}: \mathcal{X} \to
 \Delta_{1+\varepsilon} \,(=  \{\,z\in \Bbb C\,;\,|z|< 1 + \varepsilon\,\})
$$ 
of $X$ as in ~\cite{T1}, so that the fiber
$\mathcal{X}_1 := \operatorname{pr}^{-1}(1)$  is just $X$.
Then we take a  holomorphic embedding
$$
\mathcal{X} \subset \Delta \times \Bbb P^{N-1} (\Bbb C ),
$$
with $\operatorname{pr}=\pi_1$ 
such that $\pi_2^*\mathcal{O}_{\Bbb P^{N-1} (\Bbb C )}(1)$ coincides with the relative 
canonical sheaf $\mathcal{K}_{\mathcal{X}/\Delta}$ on the regular part 
$\mathcal{X}_{\operatorname{reg}}$ of 
$\mathcal{X}$,
where
$\pi_1$ (resp.~$\pi_2$) denotes the
restriction to $\mathcal{X}$ of the projection 
of the product $\Delta \times \Bbb P^{N-1} (\Bbb C )$
to the first (resp.~second) factor. Let 
$$
\psi :\; \Bbb C^* \to \operatorname{SL}(N,\Bbb C )
$$
be an algebraic group homomorphism chosen in such a way that the induced action of 
$ t \in \Bbb C^*$, $ |t|\leq 1$, 
 on $\Delta \times \Bbb P^{N-1} (\Bbb C )$ 
defined by 
$$
\Delta \times \Bbb P^N (\Bbb C ) \owns (z,p) \;\mapsto\; 
(tz, \psi (t)\cdot p)\in \Delta \times \Bbb P^{N-1} (\Bbb C )
$$
maps $\mathcal{X}$ into $\mathcal{X}$. 
Then for the Fubini-Study form $\omega_{\operatorname{FS}}$ on $\Bbb P^n$,
let $f = f(s)$ be the real-valued function defined by
$$
f (s): = \; \kappa (\psi (\exp (s))^* \omega_{\operatorname{FS}}),
\leqno{(1.1)}
$$
where $\kappa$ denotes the K-energy map. Recall that $\lim_{s\to -\infty} \dot{f}(s)$ 
is just the real part of the generalized Futaki invariant of 
the central fiber $\mathcal{X}_0$, 
and if $X$ admits a K\"ahler-Einstein metric, then $X$ is weakly K-stable 
in the sense of Tian (denoted
{\it enegy-theoretically K-stable} in this paper), i.e.,
the limit is always negative 
for all nontrivial special degenerations. 
Lemma 4.8 in this paper shows that Donaldson's K-stability for a general $L$ can be characterized 
also energy-theoretically, where in the definition of $f_m$ in (4.7), 
the K-energy appearing in the expression (1.1) for $f$ is 
 replaced by the logarithm of the Chow norm.

\medskip
This paper is organized as follows:  In Section 2, we fix notation by
defining K-stability (by Donaldson). In Section 3, by \cite{M3}, we describe the asymptotic 
behavior of the $k$-th weighted balanced metric $\omega_k$, 
 as $k\to \infty$.
In Section 4, we study the relationship between 
the Futaki invariant of a test configuration
and  the asymptotic behavior of the Chow norm of fibers (cf.~Lemma 4.8).
Finally in Section 5, based on the preceding sections, a proof for Main Theorem
will be given by developing the Chow norm method in \cite{M1} and \cite{M2}, 
where the results of 
Phong and Sturm \cite{PS1} are used to estimate the second derivative of the 
Chow norm.

\section{K-stability}

In this paper, by a {\it polarized algebraic manifold} $(M,L)$, we mean
a pair of a smooth projective algebraic variety $M$, 
defined over $\Bbb C$, and a very ample line bundle $L$
over $M$. Let $H$ be the maximal connected linear algebraic subgroup 
of the identity component $\operatorname{Aut}^0(M)$ of the
group of all holomorphic automorphisms of $M$,
so that $\operatorname{Aut}^0(M)/H$ is an 
Abelian variety (cf.~\cite{F1}).
Replacing $L$ by its suitable positive integral multiple if necessary, 
we can choose an $H$-linearization of $L$ (cf. \cite{Mu}).
Fix the natural action of the group $T := \Bbb C^*$ on the complex affine line 
 $\Bbb A^1 := \{\,z \,;\, z\in \Bbb C\,\}$ by multiplication of complex numbers,
$$
T\times \Bbb A^1 \to \Bbb A^1,
\qquad (t, z) \mapsto
tz.
$$
Let $\pi : \mathcal{M} \to \Bbb A^1$ be a $T$-equivariant projective morphism
between complex varieties with an invertible sheaf
$\mathcal{L}$ on $\mathcal{M}$,
relatively very ample over $\Bbb A^1$,
where
the algebraic group $T$ acts on $\mathcal{L}$, linearly on fibers, lifting
the $T$-action on $\mathcal{M}$. 
For each $z \in \Bbb A^1$, we put
$$
\mathcal{L}_z := {\mathcal{L}}_{|\mathcal{M}_z},
$$ 
where $\mathcal{M}_z := \pi^{-1}(z)$ denotes the scheme-theoretic fiber
of $\pi$ over $z$.
Then the following notion of a test configuration is defined by Donaldson \cite{D1}
(see special degenerations by Tian \cite{T1}).
Actually, the pair $(\mathcal{M},\mathcal{L})$ with a flat family 
$$
\pi : \mathcal{M} \to \Bbb A^1
$$ 
is called a {\it test configuration\/} for $(M,L)$, if
for some positive integer $\ell$, there exist 
the following isomorphisms of polarized algebraic manifolds
$$
(\mathcal{M}_z, {\mathcal{L}}_z )
\cong (M, \mathcal{O}^{}_M( L^{\ell})),
\qquad 0 \neq z \in \Bbb A^1.
\leqno{(2.1)}
$$

\medskip\noindent
In the special case when $\mathcal{M} = M\times \Bbb A^1$, 
a test configuration is called 
a {\it product configuration},
where for such a configuration, $T$ does not necessarily act on the first factor $M$
trivially.

\medskip
Given a test configuration $\pi : \mathcal{M} \to \Bbb A^1$ for $(M,L)$,
we consider the vector bundles  $E_m$
over $\Bbb A^1$ by
$$
\mathcal{O}^{}_{\Bbb A^1}(E_m ) \; = \; \pi_*\mathcal{L}^m,
\qquad
m =1,2,\dots,
$$
associated to the 
direct image sheaves $\pi_*\mathcal{L}^m$.
Then $E_m$ admits a natural $T$-action $\rho_m : T \times E_m \to E_m$ induced by 
the $T$-action on $\mathcal{L}$. Consider the 
fibers $(E_m)_z$, $z \in \Bbb A^1$, of the bundle
$E_m$ over $z$. Since the fiber
$$
(E_m)^{}_0 \; =\; (\pi_* \mathcal{L}^m)_0 \otimes \Bbb C
$$
over the origin is preserved by the $T$-action $\rho_m$, we can talk about the weight 
$w_m$ of the $T$-action
on $\det\, (E_m)^{}_0$. Put $n := \dim_{\Bbb C} M$, and we consider the 
degree $d_m$ of the image of the Kodaira embedding
$$
\Phi_{|L^{\ell m}|}: M \hookrightarrow \Bbb P^*(V_m),
\leqno{(2.2)}
$$
where $\Bbb P^*(V_m)$ is the set of all hyperplanes in
$V_m :=  H^0(M, \mathcal{O}^{}_M (L^{\ell m}))$
 through the origin.
Put $N_m := \dim\, (E_m)^{}_0 = \dim V_m$.
Then for 
$m \gg 1$, 
$$
\begin{cases}
\;\;\; N_m \; &= \; a_n m^n + a_{n-1}m^{n-1}+ \dots + a_1 m + a_0,\\
\;\;\; w_m \; &= \; b_{n+1} m^{n+1} + b_n m^n + \dots\, + b_1 m + b_0,
\end{cases}
\leqno{(2.3)}
$$
for some rational numbers $a_i$, $b_j \in\Bbb Q$ independent of the choice of $m$.
Note here that $a_n = \ell^n c_1(L)^n[M]/n! >0$. If $m \gg 1$, then we have
$$
\frac{w_m}{m\, N_m} \; =\; F_0 + F_1 m^{-1}+ F_2 m^{-2} + \dots ,
\leqno{(2.4)}
$$
with coeficients $F_i = F_i (\mathcal{M}, \mathcal{L} ) \in \Bbb Q$ 
independent of the choice of $m$. 
In particular $F_1 = F_1 (\mathcal{M}, \mathcal{L})$ is called the {\it Futaki invariant} 
for the test configuration. In contrast to the energy-theoretic one by Tian, the following K-stability is given by Donaldson \cite{D1}: 

\medskip\noindent
{\em Definition $2.5$.}
(i) $(M,L)$ is said to be {\it K-semistable}, 
if the inequality $F_1 (\mathcal{M}, \mathcal{L} ) \leq 0$ holds
for all test configurations 
$(\mathcal{M}, \mathcal{L})$ for $(M, L)$.

\smallskip\noindent
(ii) Let $(M,L)$ be K-semistable. Then  $(M, L)$ is said to be {\it K-stable},
if for every test configuration $(\mathcal{M}, \mathcal{L})$ for $(M, L)$,
it reduces to a product configuration if and only if $F_1(\mathcal{M}, \mathcal{L} )$ vanishes.   

\medskip
In this paper, we fix once for all a test configuration $(\mathcal{M}, \mathcal{L})$
of a polarized algebraic manifold $(M, L)$ which admits a K\"ahler metric $\omega_{\infty}$ in $c_1(L)_{\Bbb R}$ of constant 
scalar curvature. 
Obviously $F_i(\mathcal{M},\mathcal{L}^j)$ coincides with $F_i(\mathcal{M},\mathcal{L})$
for all positive integers $i$ and $j$, and hence
to discuss $K$-stability of $(M,L)$, replacing $\mathcal{L}$ by its suitable positive multiple 
if necessary, we may assume that
$\dim H^0(\mathcal{M}_z, {\mathcal{L}}^{\,m}_z)=
\dim H^0(\mathcal{M}_0, {\mathcal{L}}^{\,m}_0)$ and that the natural homomorphisms 
$$
\otimes^m H^0(\mathcal{M}_z, {\mathcal{L}}_z)
\to H^0(\mathcal{M}_z, {\mathcal{L}}^{\,m}_z),
\qquad m=1,2,\dots,
$$ 
are surjective for  all $z\in \Bbb A^1$.
We can see this easily by the fact that, if $\mathcal{L}$ is replaced by its  very high multiple 
while $L$ is fixed, then $\ell$ becomes large so that the assumptions above are
automatically satisfied (cf. \cite{Mum}; see also \cite{M4}, Remark 4.6).

\medskip
Finally, as remarked in \cite{D3}, Lemma 2, 
we have the following theorem of equivariant trivialization
for $E_m$: 

\medskip\noindent
{\bf Fact 2.6.} 
{\em  Let $H_1$ be a  Hermitian metric on the vector space $(E_m)_z$ at $z=1$.
Then there is a $T$-equivariant trivialization
$$
E_m \cong \Bbb A^1 \times (E_m)_0
\leqno{(2.7)}
$$
taking $H_1$ to a Hermitian metric, denoted by $H_0$, on the central fiber $(E_m)_0$ 
which is preserved by the action of 
$S^1 \subset \Bbb C^*\, (=T)$ on $(E_m)_0$.}
  
\section{Asymptotic behavior of weighted balanced metrics}

Now choose a Hermitian metric $h_{\infty}$ for $L$
such that $\omega_{\infty} = c_1(L, h_{\infty})$.
Let $\ell$ be as in (2.1).
Following \cite{M3}, Section 2, we here study the 
asymptotic behavior of the weighted balanced metrics for polarized 
algebraic manifolds $(M, L_{}^{m \ell })$ as
$m \to \infty$. For the linear algebraic group $H$ 
in the previous section, choose the maximal compact subgroup $K$ of $H$
such that $\omega_{\infty}$ is $K$-invariant (cf. \cite{L1}).
Then for the identity component $Z$ of the center of $K$, take its
complexfication $Z^{\Bbb C}$ in $H$.
For the $H$-linearization of $L$ in the previous section,
there exist mutually distinct characters
$\chi_{m;1}$, $\chi_{m;2}$, \dots, $\chi_{m;\nu_m}
\in \operatorname{Hom}(Z^{\Bbb C},\Bbb C^*)$ such that
the vector space $V_m$ is written as a direct sum
$$
V_m \;\, =\;\, \bigoplus_{i=1}^{\nu_m} \; V(\chi_{m;i}),
$$
where $V(\chi ):= \{\,\sigma \in V_m\,;\, g\cdot \sigma = \chi (g) \sigma 
\text{ for all $g \in Z^{\Bbb C}$}\}$ 
for all $\chi \in \operatorname{Hom}(Z^{\Bbb C},\Bbb C^*)$.
Let $\frak{z}$ be the real Lie subalgebra of $H^0(M, \mathcal{O}(T^{1,0}M))$ 
corresponding to the real Lie subgroup $Z$ of $\operatorname{Aut}(M)$.
Put $\hat{\frak{z}} := \sqrt{-1}\,\frak{z}$. 
Let $h$ be a Hermitian for $L$ such that $\omega = c_1(L; h)$ is a 
$K$-invariant K\"ahler form.
Define a $K$-invariant Hermitian pairing $\langle\; , \, \rangle^{}_h$ 
for $V_m$ by
$$
\langle \sigma, \sigma' \rangle^{}_h \; := \; \int_M (\sigma,\sigma')^{}_{h}\, \omega^{n},
\qquad \sigma, \sigma' \in V_m,
\leqno{(3.1)}
$$
where $(\sigma,\sigma')^{}_{h}$ denotes the pointwise Hermitian inner product 
of $\sigma$, $\sigma'$ by the $m$-multiple of $h$.
Then by this Hermitian pairing $\langle\; , \, \rangle^{}_h$, 
we have
$$
V(\chi_{m;i})\; \perp\; V(\chi_{m;j}),
\qquad i \neq j.
$$
Put $n_{m;i} := \dim_{\Bbb C} V(\chi_{m;i})$.
Let $P_m$ be the set of all pairs $(i, \alpha )$ of integers
such that $1 \leq i \leq \nu_m$ and $1 \leq \alpha \leq n_{m;i}$.
For the Hermitian pairing in (3.1), we say that an orthonormal basis 
$\{\, \sigma_{i,\alpha}\,;\,(i, \alpha ) \in P_m\,\}$ for $V_m$ is
{\it admissible} if $\sigma_{i,\alpha} \in V(\chi_{m;i})$ for all 
$(i,\alpha ) \in P_m$.
Fixing an admissible orthonormal basis 
$\{\, \sigma_{i,\alpha}\,;\,(i, \alpha ) \in P_m\,\}$ of 
$V_m$ with $\langle\; , \, \rangle^{}_h$, 
we now define $Z_m (\omega, \mathcal{Y}, x)$ to be
$$
(n!/m^n)\,
\Sigma_{i=1}^{\nu_m}\Sigma_{\alpha =1}^{n_{m;i}}
\exp\{-(\chi_{m;i})_*(\mathcal{Y}) + 2 x_i \}\, |\sigma_{i,\alpha}|^2_{h},
\leqno{(3.2)}
$$
for each $\mathcal{Y} \in \hat{\frak{z}}$ and 
$x =(x_1,x_2, \dots, x_{\nu_m})\in \Bbb R^{\nu_m}$,
where we put $|\sigma|_{h}^{\,2} := (\sigma,\sigma)_{h}$ for all $\sigma \in V_m$,
and $(\chi_{m;i})_*: \hat{\frak{z}} \to \Bbb R$, $i=1,2,\dots$, 
denote the differentials at $g =1$ of the restriction to $\hat{\frak{z}}$ 
of the characters $\chi_{m;i}:Z^{\Bbb C}\to\Bbb C^*$. 
Put $r_0 := n\{2 c_1(L)^n[M]\}^{-1}\{ c_1(L)^{n-1}c_1(M)[M]\}$, and consider  
$$
B_m := \{\, x = (x_1, x_2, \dots, x_{\nu_m})\in \Bbb R^{\nu_m}\, ;\,
\| x\| \leq q^2\,\},
$$
where $q := m^{-1}$ and $\|x\|:= (\Sigma_{i=1}^{\nu_m} \, n_{m;i}\,x_i^{\,2})^{1/2}$.
Then fixing a sufficiently large positive integer $k$, 
we see from \cite{M0}, Theorem B, that there exist 
vector fields $\mathcal{Y}_j\in \hat{\frak{z}}$, 
real numbers $r_j\in \Bbb R$, $j = 1,2,\dots, k$, and a $K$-invariant Hermitian metric $u_m$ for $L$ such that
\begin{align*}
Z_m (v_m, \mathcal{Y}, 0)\, &=\, 
(1+ \Sigma_{j=0}^k \,r_j q^{j+1}) \, +\, O(q^{k+2}),
\tag{3.3}\\
u_m \to h_{\infty} \text{ in}&\text{ $C^{\infty}$, \, as $m \to \infty$,}
\tag{3.4}
\end{align*}
where $v_m := c_1(L; u_m)_{\Bbb R}$ and 
$\mathcal{Y} := \Sigma_{j=1}^k \,q^{j+2} \mathcal{Y}_j$.
 In view of
the definition of 
$\delta_0$ in \cite{M2}, Step 5,
the proof of Lemma 3.4 in \cite{M1} 
allows us to make a perturbation of  $h_m$ via the action of $\exp (\frak{p}_m'' )$
(see \cite{M2} for the definition of $\frak{p}_m'' $)
to obtain a critical point for the Chow norm.
Then by (3.3) and (3.4), we obtain from \cite{M0}, pp.574--576, a $K$-invariant Hermitian metric $h_m$ for $L$
such that, for some $b_m = (b_{m;1}, b_{m;2}, \dots, b_{m;\nu_m})\in B_m$, 
\begin{align*}
Z_m (\omega_m, \mathcal{Y}, b_m&)\,\, =\,\, 
1+ \Sigma_{j=0}^k \,r_j q^{j+1},
\qquad k\gg 1,\tag{3.5}\\
h_m \to h_{\infty} \text{ in }&C^{\infty}, \, \text{ as $m \to \infty$,}
\tag{3.6}
\end{align*}
where we set $\omega_m := c_1(L; h_m)$. 
For an admissible orthonormal basis $\{\,\sigma_{i,\alpha}\,;\, (i,\alpha ) \in P_m\,\}$ 
for $V_m$ with the pairing $\langle\; , \, \rangle^{}_{h_m}$,  
by setting
$$
\beta_{m;i}\,:=\,
\exp \{ (\chi_{m;i})_*(\mathcal{Y}) + 2 b_{m;i} \}\, -1,
\leqno{(3.7)}
$$
we see from (3.5) and (3.7) the following: 
$$
(n!/m^n)\,\Sigma_{i=1}^{\nu_m}\Sigma_{\alpha = 1}^{n_{m;i}}\;
(1+\beta_{m;i})\, |\sigma_{i,\alpha}|_{h_m}^{\,2} \; = \;
 1+ \Sigma_{j=0}^k \,r_j q^{j+1},
\leqno{(3.8)}
$$
where, in view of \cite{M2}, Lemma 2.6, there exists a positive constant $C_1$ 
independent of the choice of $m\gg 1$ and $i$ such that
$$
|\beta_{m;i}| \; \leq \; C_1\, q^2
\qquad \text{for all $m \gg 1 $ and $i$.}
\leqno{(3.9)}
$$
Then by (3.8) and (3.9), we obtain
$$
\frac{\sqrt{-1}}{2\pi}\,\partial\bar{\partial}(
\Sigma_{i=1}^{\nu_m}\Sigma_{\alpha = 1}^{n_{m,i}}\, |\sigma_{i,\alpha}|_{h_m}^{\,2} )
\, -\,m\, \omega_m \; =\; O(q^2).
\leqno{(3.10)}
$$

\section{The Chow norm and the Futaki invariant}

In this section, we fix a Hermitian metric $H_1$ on $V_m$,
where 
$(E_m)_s$ at $s =1$, denoted by $(E_m)_1$, is identified with $V_m$.
By the trivialization (2.7), $H_1$ induces a Hermitian metric $H_0$
on $(E_m)_0$. Then
$$
W_m := \{\operatorname{Sym}^{d_m}((E_m)_0)\}^{\otimes n+1}
\leqno{(4.1)}
$$
admits the Chow norm (cf. Zhang \cite{Z}, 1.5; see also \S 4 in \cite{M0})
$$
W_m^* \owns\, w\;\mapsto \; \|w\|_{\operatorname{CH}(H_0 )}\, \in \Bbb R_{\geq 0}.
$$
Choose an element $\hat{M}_m$ of $W^*_m$ such that the corresponding point 
$[\hat{M}_m]$ in $\Bbb P^*(W_m)$ is the Chow point for the 
reduced effective algebraic cycle 
$$
\gamma_1 \; :=\; \Phi_{|L^{\ell m}|}(M)
$$ 
on $\Bbb P^*((E_m)_0)$.
Here each $(E_m)_s$, $s \neq 0$,  is identified with $(E_m)_0$ via
the trivialization (2.7), and by letting $s=1$, we regard $\Phi_{|L^{\ell m}|}(M)$
on $\Bbb P^*(V_m)$ 
as the algebraic cycle $\gamma_1$ on $\Bbb P^*((E_m)_0)$.
Since the $T$-action on $E_m$ 
preserves $(E_m)_0$,
we have a representation
$$
\psi_m : T \to \operatorname{GL}((E_m)_0)
\leqno{(4.2)}
$$
induced by the $T$-action on $E_m$. Note that this $T$-action on $(E_m)_0$ 
naturally induces a $T$-action on $\Bbb P^*((E_m)_0)$.
By the complete linear systems $|\mathcal{L}_s^m|$, $s \in \Bbb A^1$, 
we have the relative Kodaira embedding
$$
\mathcal{M} \;\hookrightarrow \; \Bbb P^*(E_m) ,
$$
over $\Bbb A^1$, where by (2.6) the projective bundle $\Bbb P^*(E_m)$ 
over $\Bbb A^1$
is viewed as  product bundle $\Bbb A^1 \times \Bbb P^* ((E_m)_0)$.
Then each fiber $\Bbb P^*(E_m)_s$ of $\Bbb P^*(E_m)$ over $s \in \Bbb A^1$
is naturally identified with $\Bbb P^* ((E_m)_0)$, so that
all $\mathcal{M}_z$, $z \in \Bbb A^1$, are regarded as subschemes 
of $\Bbb P^* ((E_m)_0)$. Then
$$
 \mathcal{M}_z \; =\; \psi_m (z) \cdot \mathcal{M}_1,
\qquad z \in \Bbb C^*,
\leqno{(4.3)}
$$
where on the right-hand side, the element $\psi_m (s)$ in $\operatorname{GL}((E_m)_0)$
acts naturally on $\Bbb P^* ((E_m)_0)$ as a projective linear transformation.
Note that $\mathcal{M}_1$ is nothing but $\gamma_1$ as an algebraic cycle,
and that $\mathcal{M}_0$ is preserved by the $T$-action on $\Bbb P^* ((E_m)_0)$.

\smallskip
Let us now consider the $N_m$-fold covering 
$\hat{T}:= \{\, \hat{t} \in \Bbb C^*\,\}$
of the algebraic torus $T:= \{\,t\in \Bbb C^*\,\}$ by setting
$$
t \; =\; \hat{t}^{N_m}
$$ 
for $t$ and $\hat{t}$. 
Then the mapping
$\psi_m^{\operatorname{SL}} : \hat{T} \to \operatorname{SL}((E_m)_0)$
defined by
$$
\psi_m^{\operatorname{SL}} (\hat{t}) := 
\frac{\psi_m (\hat{t}^{N_m})}{\det(\psi_m (\hat{t}))}
= \frac{\psi_m (t)}{\det(\psi_m (\hat{t}))}, 
\qquad \hat{t} \in \hat{T},
\leqno{(4.4)}
$$
is also an algebraic group homomorphism. 
Consider the quotient group 
$G_m
:=\operatorname{SL}((E_m)_0)/\Pi_m$,
 where 
$\Pi_m := \{\,\zeta^{\alpha} 
\operatorname{id}\,;\, \alpha = 1,2,\dots,N_m\,\}$
for a primitive $N_m$-th root $\zeta$ of unity.
We then define an algebraic group homomorphism $\hat{\psi}_m : T \to G_m$ by 
sending each $t\in T$ to
$$
\hat{\psi}_m (t) :\; \text{natural image of }\psi_m^{\operatorname{SL}}
 (\hat{t}) \text{ in $G_m$}.
$$
Consider
${\psi}_m (t)$, $\hat{\psi}_m (t)$, 
$\psi_m^{\operatorname{SL}} (\hat{t})$ above. Then these all
induce exactly the same projective linear transformation on $\Bbb P^*((E_m)_0)$.
Let $\gamma_t$ be the algebraic cycle on $\Bbb P^*((E_m)_0)$
obtained as the image of $\gamma_1$ by this 
projective linear transformation.
Then as $t \to 0$, we have a limit algebraic cycle
$$
\gamma_0 \; :=\; \lim_{t \to 0}\; \gamma_{t}
\leqno{(4.5)}
$$
on $\Bbb P^*((E_m)_0)$.   
To have another understanding of $\gamma_z$, 
$z \in \Bbb A^1$, recall that we  can regard each $\mathcal{M}_z$  as a subscheme 
$$
\mathcal{M}_z\; \hookrightarrow \; \Bbb P^*((E_m)_0),
\qquad z \in \Bbb A^1.
$$  
Then by (4.3), the algebraic cycle $\gamma_z$  is nothing but $\mathcal{M}_z$ viewed just  
as an algebraic cycle on $\Bbb P^*((E_m)_0)$ 
counted with multiplicities. In particular, $\gamma_0$
is the $T$-invariant algebraic cycle on $\Bbb P^*((E_m)_0)$
associated to the subscheme $\mathcal{M}_0$ counted with multiplicities. 

\smallskip
By $\hat{M}_m^{(0)}\in W_m^*$, we denote the element 
in $W_m^*$ such that the associated element $[\hat{M}_m^{(0)}]\in \Bbb P^*(W_m)$
is the Chow point for the cycle $\gamma_0$ on $\Bbb P^*((E_m)_0)$.
Then (4.5) is interpreted as
$$
\lim_{\hat{t}\to 0}\; [\,\psi_m^{\operatorname{SL}} (\hat{t} )\cdot\hat{M}_m\,]\; 
=\; [\hat{M}_m^{(0)}]
\leqno{(4.6)}
$$
in $\Bbb P^*(W_m)$. 
Here by (4.1), the group $\operatorname{GL}((E_m)_0)$ acts naturally on $W_m^*$,
and hence acts also on $\Bbb P^*(W_m)$.
We now consider the function
$$
f_m(s) \, :=\, 
\log \| \hat{\psi}_m (\exp (s))\cdot \hat{M}_m\|^{}_{\operatorname{CH}(H_0)}, 
\qquad s \in \Bbb R.
\leqno{(4.7)}
$$
Put $\dot{f}_m(s):= (df_m/ds)(s)$. The purpose of this section is to show 
the following (see Phong and Sturm \cite{PS2}, eqn 7.29, for the leading term; see also \cite{D3}, p.464--467):

\medskip\noindent
{\bf Lemma 4.8.}\,
{\em Let $a_n$ be as in $(2.3)$.
Then the function $\dot{f}_m(s)$ has a limit, as $s \to -\infty$, 
written in the following form for $m \gg 1$$:$}
\begin{align*} 
\lim_{s\to -\infty} \dot{f}_m(s) \; &=\;   (n+1)!\, a_n  (F_1 m^n + F_2 m^{n-1} + F_3 m^{n-2} + \dots )  \tag{4.9}
\\
& =\; (n+1)! \, a_n \left (\frac{w_m}{mN_m} - F_0\right ) m^{n+1}.
\end{align*}
{\em Proof$:$} \;  
Since $\gamma_0$ is preserved by the $\hat{T}$-action on $(E_m)_0$, the Chow point 
$[\hat{M}^{(0)}]$ for $\gamma_0$ is fixed by the $\hat{T}$-action on $\Bbb P^*(W_m)$,
i.e., for some $q_m\in \Bbb Z$, 
$$
\psi_m^{\operatorname{SL}} (\hat{t}) \cdot \hat{M}_m^{(0)}\; 
= \; \hat{t}_{}^{q_m}\hat{M}_m^{(0)},
\qquad t \in \Bbb C^*,
$$
where the left-hand side is $\hat{\psi}_m (t) \cdot \hat{M}^{(0)}$ 
modulo the action of $\Pi_m$. Since the $\hat{T}$-action on $W_m^*$ is diagonalizable, 
we can write $\hat{M}_m$ in the form
$$
\hat{M}_m \; =\; \Sigma_{\nu =1}^N \; w_{\nu},
\leqno{(4.10)}
$$
where  $\,0 \neq w_{\nu}\in W_m^*$, $\nu =1,2,\dots, N$, 
are such that,  for an increasing sequence of integers 
$e_1 <e_2 < \dots < e_N$, the equality
$$
\psi_m^{\operatorname{SL}} (\hat{t})\cdot w_{\nu} 
= \hat{t}^{e_{\nu}}w_{\nu}
\leqno{(4.11)}
$$ 
holds for all $\nu\in \{1,2,\dots,N\}$ and $\hat{t}\in \hat{T}$.
In particular, in view of (4.6), we can find a complex number $c \neq 0$ 
such that
$$
\hat{M}_m^{(0)}\; =\; c\,w_1,
$$
and hence $q_m$ coincides with $e_1$.
Then by (4.10) and (4.11), it is easy to check that
$$
\lim_{s\to -\infty} \dot{f}_m(s) \; \left (=\; \frac{e_1}{N_m} \right )\;  =\;
\frac{q_m}{N_m}.
\leqno{(4.12)}
$$
Hence it suffices to show that $q_m/N_m$ admits the asymptotic expansion 
as in the right-hand side of (4.9) above. Consider the graded algebra
$$
\bigoplus_{k=0}^{\infty}\; (E_{km})_0,
$$
where  via $\psi_m^{\operatorname{SL}}$, 
the group  $\hat{T}$ acts on $(E_{m})_0$ 
and hence on $(E_{km})_0$.
Then by \cite{Mu0}, Proposition 2.11, 
the weight $p_{k}$ for the $\hat{T}$-action on $\det (E_{km})_0$ 
satisfies the following:
$$
p_k \, +\, \frac{q_m}{(n+1)!}\, k^{n+1}\; = \; O(k^n), \qquad k \gg 1,
\leqno{(4.13)}
$$
i.e., there exists a constant $C >0$
independent of $k$, possibly depending on $m$, 
such that the left-hand side of (4.13) 
has absolute value bounded by $C k^n$ for  positive integers $k$.
Recall the definition of $w_{km}$ and $w_m$ in Section 2. 
Then by the expression of 
$\psi_m^{\operatorname{SL}}$ in (4.4), the weight $p_k$ for $\det (E_{km})_0$ 
induced by the $\hat{T}$-action on $(E_m)_0$ via $\psi_m^{\operatorname{SL}}$
is expressible as
$$
p_k \; =\; N_m w_{km}\; -\; k\, w_m N_{km}.
\leqno{(4.14)}
$$
Here the term $N_m w_{mk}$ in the right-hand side of (4.14) 
is the weight in $\hat{t}$ for $\det (E_{km})_0$
induced from the action on $(E_m)_0$ by the numerator $\psi_m(t)$ 
in (4.4), since it
is nothing but the weight in $\hat{t}$ for the 
action of $\psi_{mk}(t)$ on 
$\det (E_{km})_0$, 
while in view of the natural surjective homomorphism
$$
S^k ((E_m)_0) \to (E_{km})_0,
$$
the term $k\, w_m N_{km}$  is just
the weight in $\hat{t}$ induced from the scalar action on $(E_m)_0$ by
the denominator of (4.4).
Then for $k \gg 1$, by (4.14) and (2.4), we obtain
\begin{align*}
p_k \; &=\; (km)\, N_m N_{km} 
\left \{\, \frac{w_{km}}{(k m) N_{km}}\, -\, \frac{w_m}{m N_m}\,\right \}\\
&=\; -\,(km)\, N_m N_{km}\{\, (F_1m^{-1} + F_2 m^{-2} + F_3 m^{-3} + \dots )\, +\,O(k^{-1}) \,\}\\
&=\; -\, k^{n+1}a_n N_m \{\, (F_1 m^n + F_2 m^{n-1} + F_3 m^{n-2} + \dots )\, 
+\,O(k^{-1}) \,\},
\end{align*}
where the last equality ifollows from (2.3) applied to 
$km$.
 Then by comparing this with (4.13),  and then by (2.4), we obtain
\begin{align*}
\frac{q_m}{N_m}\; &=\;  (n+1)!\, a_n  (F_1 m^n + F_2 m^{n-1} + F_3 m^{n-2} + \dots )\\
&= \; (n+1)! \, a_n \left (\frac{w_m}{mN_m} - F_0\right ) m^{n+1}.
\end{align*}
as required.
\qed

\section{Proof of Main Theorem}

In this section, by using the notation in (3.1), we choose $\langle\; , \, \rangle^{}_{h_m}$ 
as the Hermitian metric $H_1$ for $V_m$ in Section 4,
where $h_m$ is as in (3.6).
For the corresponding Chow norm
$$
W_m^* \owns\, w\;\mapsto \; \|w\|_{\operatorname{CH}(H_0 )}\, \in \Bbb R_{\geq 0},
$$
we consider the real-valued function $f_m$ on $\Bbb R$ as in (4.7).
For the one-parameter group $\hat{\psi}_m : \hat{T} \to G_m$, 
the vector space $(E_m)_0$ admits an orthonormal basis $\mathcal{T}:=\{\, \tau_1, \tau_2, \dots, \tau_{N_m}\}$ such that, for some $e_i\in \Bbb Q$ 
with $\Sigma_{\alpha =1}^{N_m} e_{\alpha} = 0$, we have
$$
\hat{\psi}_m (t)\, \tau_{\alpha} \; \equiv\; t^{e_{\alpha}} \tau_{\alpha},
\qquad t \in T,
\leqno{(5.1)}
$$
modulo the action of $\Pi_m$.
By the associated Kodaira embedding
$\mathcal{M}_0 \hookrightarrow \Bbb P^*((E_m)_0 )$, we regard $\mathcal{M}_0$ 
as a subscheme 
$$
\mathcal{M}_0 \; \hookrightarrow \; \Bbb P^{N_m-1}(\Bbb C ),
\qquad \, p\, \mapsto (\tau_1(p):\tau_2(p):\dots :\tau_{N_m}(p)), 
$$
where we identify $\Bbb P^*((E_m)_0 )$ 
with 
$$
\Bbb P^{N_m-1}(\Bbb C ):= \{(z^{}_1:z^{}_2:\dots :z^{}_{N_m})\}
$$ 
by the basis $\mathcal{T}$. Since we regard $(E_m)_1$ just as  $V_m$, the identification (2.7)
allows us to obtain a basis
$\mathcal{T}'\; :=\; \{\, \tau'_1, \tau'_2, \dots, \tau'_{N_m}\}$ 
for $V_m$  
corresponding to the basis $\mathcal{T}$ for $(E_m)_0$.
Note that this basis $\mathcal{T}'$ is orthonormal
with respect to the Hermitian metric $H_1 = \langle\; , \, \rangle^{}_{h_m}$. 
Then the Kodaira embedding 
$\Phi_{|L^{\ell m}|}: M\,  (=\,\mathcal{M}_1)\hookrightarrow \Bbb P^* (V_m)$
is given by
$$
M \; \hookrightarrow \; \Bbb P^{N_m-1}(\Bbb C ),
\qquad \, p\, \mapsto (\tau'_1(p):\tau'_2(p):\dots :\tau'_{N_m}(p)),
\leqno{(5.2)}
$$
where we identify 
$\Bbb P^*(V_m) = \Bbb P^{N_m-1}(\Bbb C ) = \Bbb P^*((E_m)_0)$ 
by the bases $\mathcal{T}'$ and $\mathcal{T}$.
Then the Fubini-Study form 
$\omega_{\operatorname{FS}}$ 
on $\Bbb P^{N_m-1} (\Bbb C )\, (= \Bbb P^*(V_m))$ is 
$$
\omega_{\operatorname{FS}}\; :=\; 
\frac{\sqrt{-1}}{2\pi} \,\partial\bar{\partial}\log (\Sigma_{\alpha =1}^{N_m}
|z_{\alpha}|^2).
$$ 
By (3.10), we here observe that 
$$
\omega_{\operatorname{FS}} - m \omega_m = O(q^2),
\leqno{(5.3)}
$$
on $M$. For the function $f_m(s)$ in (4.7), we first give an estimate 
of the fist derivative $\dot{f}_m(0)$.
In view of \cite{Z} (see also \cite{M0}), 
$$
\dot{f}_m(0)\; =\; (n+1)\int_M 
\frac{\Sigma_{\alpha =1}^{N_m}\, e_{\alpha}\,
|\tau'_{\alpha}|^2_{h_m}}
{\Sigma_{\alpha =1}^{N_m}\, 
|\tau'_{\alpha}|^2_{h_m}}\,
\omega_{\operatorname{FS}}^n,
\leqno{(5.4)}
$$
where by (3.8) and (3.9), we observe that
\begin{align*}
&\Sigma_{\alpha =1}^{N_m}\, 
|\tau'_{\alpha}|^2_{h_m}\; =\; \,\Sigma_{i=1}^{\nu_m}\Sigma_{\alpha = 1}^{n_{m,i}}\;
|\sigma_{i,\alpha}|_{h_m}^{\,2} \; 
\tag{5.5} \\
&= \;
(m^n/n!)\, ( 1+ \Sigma_{j=0}^k \,r_j q^{j+1})\,
-\, \Sigma_{i=1}^{\nu_m}\Sigma_{\alpha = 1}^{n_{m,i}}\;
\beta_{m;i}\, |\sigma_{i,\alpha}|_{h_m}^{\,2}\\
&=\; (m^n/n!)\, ( 1+ \Sigma_{j=0}^k \,r_j q^{j+1})\, \{ 1 + O(q^2)\}.
\end{align*}
Now, we can rewrite (5.4) in the form
\begin{align*}
\dot{f}_m(0)\; &=\; (n+1)!\int_M
\frac{\Sigma_{\alpha =1}^{N_m}\, e_{\alpha}\,
|\tau'_{\alpha}|^2_{h_m}}{1+ \Sigma_{j=0}^k \,r_j q^{j+1}}
\,\{ 1 + O(q^2)\}
\, \omega_m^{\,n}\tag{5.6}  \\
&=\;\int_M O(q^2)
(\Sigma_{\alpha =1}^{N_m}\, e_{\alpha}\,
|\tau'_{\alpha}|^2_{h_m})\,\omega_m^{\,n},
\end{align*}
where the last equality follows from
$$
\int_M \Sigma_{\alpha =1}^{N_m}\, e_{\alpha}
|\tau'_{\alpha}|^2_{h_m}\omega_m^{\,n}\, =\,\Sigma_{\alpha =1}^{N_m}\, e_{\alpha}
\, =\,0.
$$
All weights $e_{\alpha}$ for $\hat{\psi}_m$ have absolute value
bounded by $C_2 m$ 
for some constant $C_2 > 0$ independent of both $\alpha$ and $m$, i.e.,
$$
|e_{\alpha}| \; \leq \; C_2 m,
\qquad \alpha = 1,2,\dots, N_m.
\leqno{(5.7)}
$$
In view of (5.5), $\Sigma_{\alpha =1}^{N_m}\, 
|\tau'_{\alpha}|^2_{h_m}\, = O(m^n)$.
Then by (5.6) and (5.7), 
$$
\dot{f}_m (0) \; =\; O(m^{n-1}).
\leqno{(5.8)}
$$
By \cite{Z} (see also \cite{M0}, 4.5),  $\ddot{f}_m(s) \geq 0$ for all $s\in \Bbb R$.
In (4.9), let $m \to \infty$. Then by (5.8),  we obtain $F_1 \leq 0$,
i.e., K-semistability of $(M,L)$ follows.

\medskip 
To show K-stability, we now assume that $F_1 =0$ for the test configuration 
above. Then by Lemma 4.8,
$$
\lim_{s\to -\infty} \dot{f}_m (s) = O(m^{n-1}).
\leqno{(5.9)}
$$
Here we consider the second derivative $\ddot{f}_m(s)$.
From now on, by setting $\delta := C_3 (\log m ) q$, 
we require the real number $s$ to satisfy 
$$
|s| \; \leq \; \delta, 
\leqno{(5.10)}
$$
where $C_3$ is a positive real number independent of the choice of $m$.
For local one-parameter group 
$$
\mu_{m,s} := \hat{\psi}_m(\exp (s)) \in G_m,
\qquad -\delta \leq s \leq \delta,
$$
we regard each $\mu_{m,s}$ as a linear isomorphism 
of $\Bbb C^{N_m -1} \,(=V_m)$, modulo the action by 
$\Pi_m$,
via the identification of $(E_m)_0$ with $(E_m)_1 \, (= V_m)$.
Note also that $G_m$ acts on $\Bbb P^*((E_m)_0)\, (=\Bbb P^*(V_m))$ via the 
projection
$$
\pi_{m} : \, G_m\,(= \operatorname{SL}
(N_m; \Bbb C ) /\Pi_m )\, \to\, 
\operatorname{PGL}((E_m)_0)\, (=\operatorname{PGL}(N_m;\Bbb C )).
$$
Now by Appendix, the family of K\"ahler manifolds 
$$
(M,\,q\, (\mu_{m,s}^*{\omega_{\operatorname{FS}}})_{|M}),
\qquad -\delta \leq s \leq \delta,\; m=1,2,\dots,
\leqno{(5.11)}
$$
has bounded geometry.
Let us now consider the holomorphic vector field
$\mathcal{V}^m:=(\pi_{m}\circ\hat{\psi}_m)_* (\partial/\partial s)$ on $\Bbb P^{N_m-1}(\Bbb C)$ 
which generates the local one-parameter group 
$\pi_m (\mu_{m,s} )$, $-\delta \leq s \leq \delta$.
For each $s$, consider the holomorphic tangent bundle $TM_{s}$
of $M_{s} := \mu_{m,s} (M)$,
where $M$ is viewed as a subvariety of $\Bbb P^{N_m-1}(\Bbb C)$ by (5.2).
Metrically, for the orthogonal complement $TM_{s}^{\perp}$ 
of $TM_{s}$ in $T\Bbb P^{N_m-1}(\Bbb C)_{|M_{s}}$ 
by the metric $\omega_{\operatorname{FS}}$, 
we can regard
the normal bundle of $M_{s}$ in $\Bbb P^{N_m-1}(\Bbb C)$ 
as the subbundle $TM_{s}^{\perp}$ 
of $T\Bbb P^{N_m-1}(\Bbb C)_{|M_{s}}$. 
 Hence $T\Bbb P^{N_m-1}(\Bbb C)_{|M_{s}}$ is differentiably a
direct sum $TM_{s}\oplus TM_{s}^{\perp}$, and we can uniquely write 
$$
{\mathcal{V}^m}_{|M_{s}} \; = \; \mathcal{V}^m_{TM_{s}} \,+\, 
\mathcal{V}^m_{TM_{s}^{\perp}},
\leqno{(5.12)}
$$ 
where $\mathcal{V}_{TM}$ and $\mathcal{V}_{TM^{\perp}}$ are smooth sections of 
$TM_{s}$ and $TM_{s}^{\perp}$, respectively.
Consider the exact sequence of holomorphic vector bundles
$$
0 \rightarrow TM_{s} \rightarrow T\Bbb P^{N_m-1}(\Bbb C)_{|M_{s}} 
\rightarrow TM_{s}^{\perp} \rightarrow 0
$$ 
over $M_{s}$.
Then the pointwise estimate (cf. \cite{PS1}, (5.16)) for the second fundamental form for this
exact sequence is valid also in our case (cf. \cite{M1}, Step 2), and as in \cite{PS1}, (5.15), we obtain the
inequality
$$
 \int_{M_{s}} |\mathcal{V}^m_{TM_{s}^{\perp}}|^2_{\omega_{\operatorname{FS}}}\, \omega_{\operatorname{FS}}^n
\; \geq \; C_4\,
\int_{M_{s}} |\bar{\partial}\mathcal{V}^m_{TM_{s}^{\perp}}|^2_{\omega_{\operatorname{FS}}}
\,\omega_{\operatorname{FS}}^n,
\leqno{(5.13)}
$$
where $C_4$ is a positive real constant independent of the choice of $m$.
The second derivative $\ddot{f}_m (s)$ is (see for instance \cite{M0}, \cite{PS1}) 
given by
$$
\ddot{f}_m (s) \; =\;\int_{M_{s}} \,
  | \mathcal{V}^m_{TM_{s}^{\perp}} |^2_{\omega_{\operatorname{FS}}}\,
\omega_{\operatorname{FS}}^n\; \geq \; 0.
\leqno{(5.14)}
$$
Put $\varphi_m := (\Sigma_{\alpha = 1}^{N_m}\, e_{\alpha} |z_{\alpha} |^2)
/(m \Sigma_{\alpha = 1}^{N_m} |z_{\alpha} |^2)$ on $\Bbb P^{N_m -1}(\Bbb C )$.
Then by (5.7) $\varphi_m$, $m = 1,2,\dots$, are uniformly bounded  satisfying (cf. \cite{M1}, (4.5))
$$
i_{\mathcal{V}^m}(q\, \omega_{\operatorname{FS}})
\; =\; \frac{\sqrt{-1}}{2\pi}\, \bar{\partial}\varphi_m.
\leqno{(5.15)}
$$
For $M$ as a submanifold of 
$\Bbb P^{N_m-1}(\Bbb C )$ in (5.2), 
we consider the sheaves $\mathcal{A}^{0,j}(TM)$, 
$j=0,1,\dots$, of germs of smooth 
$(0,j)$-forms on $M$ with values in the holomorphic tangent bundle $TM$ 
of $M$, and endow $M$ with the K\"ahler metric ${q\mu_{m,s}^* \omega_{\operatorname{FS}}}_{|M}$ 
for $s$ as in (5.10).
Now on $\mathcal{A}^{(0,0)}(TM)$, we consider the operator 
$$
\Delta_{TM,s} \; :=\; - \,\bar{\partial}^{\#}\bar{\partial},
$$
where $\bar{\partial}^{\#}$ is the 
formal adjoint of $\bar{\partial}:\mathcal{A}^{(0,0)}(TM)\to
\mathcal{A}^{(0,1)}(TM)$.  
For $\Gamma :=H^0 (M,\mathcal{A}^{(0,0)}(TM))$, we consider the Hermitian $L^2$-pairing 
$$
\langle V_1, V_2\rangle^{}_{s}
 \; := \; \int_M (V_1, V_2)_{q\mu_{m,s}^* \omega_{\operatorname{FS}}}
\, (q\mu_{m,s}^* \omega_{\operatorname{FS}})^n,
\qquad V_1, V_2 \in \Gamma,
$$
where $(V_1, V_2)_{q\mu_{m,s}^* \omega_{\operatorname{FS}}}$ is the 
pointwise Hermitian pairing of $V_1$ and $V_2$ by the K\"ahler metric 
$q\mu_{m,s}^* \omega_{\operatorname{FS}}$. 
For the subspace $\frak{g}:= H^0 (M, \mathcal{O}(TM))$ of
$\Gamma$, we consider its orthogonal complement $\frak g_{s}^{\perp}$ in $\Gamma$ by
the pairing $\langle \;, \;\rangle^{}_{s}$.
Then $\mathcal{V}^m_{TM_{s}}$ in (5.12) is expressible as
$$
\mathcal{V}^m_{TM_{s}}\;=\; \mathcal{V}_{m,s}^{\circ}\,+\, \mathcal{V}_{m,s}^{\bullet},
$$
where $\mathcal{V}_{m,s}^{\circ}$ and $\mathcal{V}_{m,s}^{\bullet}$ belong to $(\mu_{m,s})_*\frak g$ and 
$(\mu_{m,s})_*\frak g_{s}^{\perp}$, respectively.
Since the left-hand side of (5.12) is holomorphic, 
$$
\bar{\partial}\mathcal{V}^m_{TM_{s}^{\perp}}\; =\; -\,\bar{\partial}\mathcal{V}^m_{TM_{s}}
\; =\;-\,\bar{\partial}\mathcal{V}_{m,s}^{\bullet}. 
\leqno{(5.16)}
$$
Since the family (5.11) has bounded geometry, the first positive eigenvalue $\lambda_1$ 
of the operator $\Delta_{TM}$ on $\mathcal{A}^{(0,0)}(TM)$ is bounded 
from below by some positive constant $C_5$ independent of the choice of $m$. Hence
$$
\int_{M_{s}}|\bar{\partial}\mathcal{V}_{m,s}^{\bullet}|^2_{q\omega_{\operatorname{FS}}}
(q\omega_{\operatorname{FS}})^n\; \geq \; 
C_5\int_{M_{s}} |\mathcal{V}_{m,s}^{\bullet}|^2_{q\omega_{\operatorname{FS}}}(q\omega_{\operatorname{FS}})^n.
\leqno{(5.17)}
$$
From (5.13), (5.16) and (5.17), it now follows that
$$
\int_{M_{s}} |\mathcal{V}^m_{TM_{s}^{\perp}}|^2_{\omega_{\operatorname{FS}}} 
\,\omega_{\operatorname{FS}}^n\; \geq \; C_4 C_5\, q \, 
\int_{M_{s}} |\mathcal{V}_{m,s}^{\bullet}|^2_{\omega_{\operatorname{FS}}}\,
\omega_{\operatorname{FS}}^n.
\leqno{(5.18)}
$$
In view of (5.14), 
$\dot{f}_m(0) - \dot{f}_m(-\delta )  = \int_{-\delta}^0\,\ddot{f}_m(s)\, ds \,\geq \,0$,
and it follows from (5.8) and (5.9) that 
\begin{align*}
O(m^{n-1}) \; &=\; \dot{f}_m (0) - \lim_{s\to -\infty}  \dot{f}_m (s)  \;
\geq \; 
\dot{f}_m (0) - \dot{f}_m (-\delta )  \\
&=\; \int_{-\delta}^0\,\ddot{f}_m(s)\, ds\; \geq \;  \ddot{f}_m(s_m) \, \delta
 \end{align*}
where $s_m$, $m =1,2,\dots$, are real numbers at which the functions
$\ddot{f}_m(s)$,  $-\delta \leq s \leq 0$, attain their minima, i.e., 
$$
\ddot{f}_m (s_m)\; = \; \min_{-\delta \leq s \leq 0} \ddot{f}_m (s).
$$
Therefore, in view of (5.14) and $\delta \,=\,O(q\log m )$, we obtain
$$
\int_{M^{}_{s^{}_m}}  | \mathcal{V}^m_{TM_{s_m}^{\perp}} |^2_{q\omega_{\operatorname{FS}}}\,
(q\omega_{\operatorname{FS}})^n\; 
 =\; 
O\left (\frac{q}{\log m}\right ),
\leqno{(5.19)}
$$
since the left-hand side is $\ddot{f}_m (s_m )/ m^{n+1}$.
Then by (5.18), 
$$
\int_{M^{}_{s^{}_m}} 
 |\mathcal{V}^{\bullet}_{m,{s^{}_m}}
 |^2_{q\omega_{\operatorname{FS}}}\,
(q\omega_{\operatorname{FS}})^n
\, =\; O\left (\frac{1}{\log m}\right ).
\leqno{(5.20)}
$$
Infinitesimally, (5.1) is written as
$\mathcal{V}^m\cdot \tau'_{\alpha} \, =\, e_{\alpha} \tau'_{\alpha}$,
$\alpha = 1,2,\dots, N_m$,
and for $\varphi_m$  in (5.15), by setting
$$
\bar{\varphi}_m :=\;  \left\{\int_{M_{s_m}}(q \omega_{\operatorname{FS}})^n
\right \}^{-1}
\int_{M_{s_m}}\varphi_m (q \omega_{\operatorname{FS}})^n,
$$
we can write the pullback
$\mu_{m, s_m}^*(\varphi_m - \bar{\varphi}_m)$  as
$$
\eta_m \; :=\; (\mu_{m,s_m}^*\varphi^{}_m)^{}_{|M}-\bar{\varphi}_m \; =\; 
\frac{\Sigma_{\alpha = 1}^{N_m}\, e_{\alpha}\, \tau'_{\alpha}\bar{\tau}'_{\alpha}
\exp (2s_m e_{\alpha})}
{m \Sigma_{\alpha = 1}^{N_m}\, \tau'_{\alpha}\bar{\tau}'_{\alpha}
\exp (2s_m e_{\alpha}) } - \bar{\varphi}_m
\leqno{(5.21)}
$$
when restricted to $M \subset \Bbb P^*(V_m)$. 
Note that $\bar{\varphi}_m$, $m = 1,2,\dots$, is a bounded sequence 
of real numbers.
Then 
we can write the uniformly bounded real-valued functions $\eta_m$ on $M$ as
$$
\eta_m \; := \frac{\Sigma_{\alpha = 1}^{N_m}\, e_{\alpha,m}\, \tau'_{\alpha}\bar{\tau}'_{\alpha}
\exp (2s_m e_{\alpha ,m})}
{m \Sigma_{\alpha = 1}^{N_m}\, \tau'_{\alpha}\bar{\tau}'_{\alpha}
\exp (2s_m e_{\alpha ,m}) },
\qquad m= 1,2,\dots,
$$
where $e_{\alpha,m}:= e_{\alpha} - m \bar{\varphi}_m $.
Put $\omega_m := 
{q\mu_{s_m}^*\omega_{\operatorname{FS}}}_{|M}$.
Hereafter, replace the sequence $s_m$, $m \gg 1$, by its suitable subsequence 
$s_{m_j}$, $j=1,2,\dots$, if necessary.  
We write $m_j$, $m_j^{-1}$, $N_{m_j}$, $s_{m_j}$, $\omega_{m_j}$, $\eta_{m_j}$, 
$e_{\alpha, m_j}$
as $m(j)$, $q(j)$, $N(j)$, $s(j)$, $\omega (j)$, $\eta (j)$, $e_{\alpha}(j)$,
respectively.
Then by Appendix, we may assume that $\omega (j)$
converges to $\omega_{\infty}$ in $C^{\infty}$ as $j \to \infty$.
Moreover, we set
$$
\begin{cases}
&\mathcal{V} (j):= \mathcal{V}^{m(j)} = (\mu_j^{-1})_* \mathcal{V}^{m(j)}, \;\;
\mathcal{V}_{TM}(j):=(\mu_j^{-1})_*\mathcal{V}^{m(j)}_{TM_{s(j)}},\\
&\mathcal{V}^{\circ}(j):= (\mu_j^{-1})_*\mathcal{V}^{\circ}_{m(j),{s(j)}},\; 
\mathcal{V}^{\bullet}(j):= (\mu_j^{-1})_*\mathcal{V}^{\bullet}_{m(j),{s(j)}},
\end{cases}
$$
where $\mu_j := \mu_{m(j), s(j)}$.
Then the  following cases 1 and 2 are possible: 

\medskip\noindent
Case 1: $I^{\circ}_j := \int_{M} 
 |\mathcal{V}^{\circ}(j)
 |^2_{\omega (j)}\,
\omega (j)^n$, $j =1,2,\dots$, 
are bounded. 
In this case, by $ |\mathcal{V}^2_{TM}(j)|^{\,2}_{\omega (j)}
=  |\mathcal{V}^{\circ}(j) |^{\,2}_{\omega (j)}
+ |\mathcal{V}^{\bullet}(j) |^{\,2}_{\omega (j)}$,
this boundedness together with (5.20)  
implies that 
$$
\int_M 
 |\mathcal{V}_{TM}(j)
 |^2_{\omega (j)}\,
\omega (j)^n,\; j=1,2,\dots, \text{ are bounded.}
\leqno{(5.22)}
$$
Note that $\omega (j) \to \omega_{\infty}$ in $C^{\infty}$, as $j \to \infty$.
Hence in view of (5.22), since $ |\mathcal{V}_{TM}(j)
 |^2_{\omega (j)} = |\bar{\partial}{\eta}(j)|^2_{\omega (j)}$ by (5.15), 
the sequence of integrals
$\int_M |\bar{\partial}{\eta}(j)|^2_{\omega_{\infty}}\omega_{\infty}^n$,
$j = 1,2,\dots$, is bounded, so that $\eta (j)$, $j =1,2,\dots$, is 
a bounded sequence in the Sobolev space $L^{1,2}(M, \omega_{\infty})$.
Now by \cite{RT}, we may assume that $n \geq 2$.
Then replacing $\eta (j)$, $j =1,2,\dots$, by its subsequence if necessary,
we may further assume the convergence
$$
\eta (j) \to \eta_{\infty} \text{ in $L^2(M, \omega^n_{\infty})$}, \qquad
\text{as $j \to \infty$,}
\leqno{(5.23)}
$$
where $\eta_{\infty}$ is a real-valued function  
in $L^2(M, \omega_{\infty})$.
Recall that the Lichnerowich operator $\Lambda_j : C^{\infty}(M)_{\Bbb C} \to C^{\infty}(M)_{\Bbb C}$ 
for the K\"ahler manifold $(M, \omega (j))$ is an 
elliptic operator, of order 4, with kernel consisting of 
all Hamiltonian functions for
the holomorphic Hamiltonian vector fields on $M$. 
Now, to  each smooth function $f \in C^{\infty}(M)_{\Bbb C}$, 
we associate a complex vector field $\frak S_{f, j}$ 
of type $(1,0)$ on $M$
such that
$$
i(\frak S_{f,j})\, \omega (j) \, = \, \frac{\sqrt{-1}}{2\pi}\,\bar{\partial}f.
$$
Note that $\frak S_{\eta (j),j} = \mathcal{V}_{TM}(j)$ by (5.15) and (5.21).
Then for the formal adjoint $\Lambda_j^{\#}: 
C^{\infty}(M)_{\Bbb C} \to C^{\infty}(M)_{\Bbb C}$ of 
the operator $\Lambda_j$, we have
\begin{align*}
&\int_M \eta (j )\,\{\Lambda_j^{\#}f\}\, \omega (j)^n\; =\;
\int_M \{ \Lambda_j \eta (j) \}\,f\,\omega (j)^n \; = \;
\langle \,\bar{\partial}\frak S_{\eta (j), j}, \bar{\partial}\frak S_{f, j} \rangle_{s(j)}\\
&=\;\langle\, \bar{\partial}\{\mathcal{V}_{TM}(j)\}, 
\bar{\partial}\frak S_{f, j} \rangle_{s(j)}
\;=\; \langle\, \bar{\partial}\{\mathcal{V}^{\bullet}(j)\}, 
\bar{\partial}\frak S_{f, j} \rangle_{s(j)},
\end{align*}
for all $f \in C^{\infty}(M)_{\Bbb C}$.
Here the last equality follows from
the identities $\mathcal{V}_{TM}(j) = \mathcal{V}^{\circ}(j) + \mathcal{V}^{\bullet}(j)$
and $\bar{\partial}\mathcal{V}^{\circ}(j) =0$. Hence, for each 
fixed $f$ in $C^{\infty}(M)_{\Bbb C}$, we obtain
$$
\begin{cases}
&\;\left | \int_M \eta (j) \,\{\Lambda_j^{\#}f\}\, \omega (j)^n\right |\; =\;
\;
\left | \langle\, \mathcal{V}^{\bullet}(j), 
\Delta_j \frak S_{f, j}  \rangle_{s(j)} \right | \\ 
&\; \leq\;  \left \{ \int_M |\Delta_j \frak S_{f, j}|^2_{\omega (j)}\, 
\omega (j)^n  \right \}^{1/2}  I^{\bullet}_j,
\end{cases}
\leqno{(5.24)}
$$
where $I^{\bullet}_j := \{\int_{M} 
 |\mathcal{V}^{\bullet}(j)
 |^2_{\omega (j)}\,
\omega (j)^n\}^{1/2}$ and
$\Delta_j :=\Delta_{TM, s(j)}$. 
Let $j \to \infty$ in (5.24). Since $I^{\bullet}_j \to 0$ by (5.20), 
and since $\omega (j ) \to \omega_{\infty}$ in $C^{\infty}$,
by passing to the limit, we see from  (5.23) and (5.24) that
$$
\int_M \eta_{\infty} \,\{\Lambda_{\infty}^{\#}f\}\, \omega_{\infty}^n
\; =\; 0
\qquad \text{for all $f \in C^{\infty}(M)_{\Bbb C}$},
$$ 
where $\Lambda_{\infty}: C^{\infty}(M)_{\Bbb C} \to  C^{\infty}(M)_{\Bbb C}$
is the Lichnerowich operator for the K\"ahler manifold $(M, \omega_{\infty})$,
and $\Lambda_{\infty}^{\#}$ is its formal adjoint.
Since $\Lambda_{\infty}$ is elliptic, any weak solution $\eta = \eta_{\infty}$ for the equation
$$
\Lambda_{\infty} \eta =  0
$$
is always a strong solution. In particular $\eta_{\infty}$ is 
a real-valed smooth function on $M$ such that the complex vector field 
$W$ of type $(1,0)$ on $M$ defined by
$i(W) \omega_{\infty} \; =\; \bar{\partial}\eta_{\infty}$
is holomorphic. 
Then by \cite{M5}, the test configuration
$\pi : \mathcal{M} \to \Bbb A^1$ is a product configuration.

\medskip\noindent
Case 2: $I^{\circ}_j \to +\infty$ as $j \to \infty$. In this case, we put 
$\hat{\mathcal{V}}_{TM}(j):= \mathcal{V}_{TM}(j)/\sqrt{I^{\circ}_j}$, 
$\hat{\mathcal{V}}^{\circ}(j):= \mathcal{V}^{\circ}(j)/\sqrt{I^{\circ}_j}$, and 
$\hat{\mathcal{V}}^{\bullet}(j):= \mathcal{V}^{\bullet}(j)/\sqrt{I^{\circ}_j}$.
Then
$$
\int_M  |\hat{\mathcal{V}}^{\circ}(j) |^2_{\omega (j)}
\,\omega (j)^n \; = \; 1,
\qquad j=1,2,\dots,
\leqno{(5.25)}
$$
where by setting
$\hat{\eta}(j): = \eta (j)/\sqrt{I^{\circ}_j}$, we see 
from $\frak S_{\eta (j),j} = \mathcal{V}_{TM}(j)$ that the complex
vector field $\hat{\mathcal{V}}_{TM}(j)$ of type $(1,0)$ on $M$ satisfies
$$
 i(\hat{\mathcal{V}}_{TM}(j))\,\omega (j) \; =\; \sqrt{-1}\,\bar{\partial}(\hat{\eta}(j)).
 \leqno{(5.26)}
$$
Since the functions $\eta_m$, $m \gg 1$, are uniformly bounded on $M$, 
and since  $\omega (j)$ converges to $\omega_{\infty}$ as $j \to \infty$, 
we obtain the convergence
$$
\hat{\eta}(j) \to 0 \text{ in $C^0(M)$}, 
\qquad\text{as $j \to \infty$.}
\leqno{(5.27)}
$$
By (5.25), replacing $\hat{\mathcal{V}}^{\circ}(j)$, 
$j = 1,2,\dots$, by its subsequence if necessary, we may assume that 
$$
\hat{\mathcal{V}}^{\circ}(j) \to \hat{\mathcal{V}}^{\circ}_{\infty}\text{ in $\frak g$},
\qquad \text{as }j \to \infty,
\leqno{(5.28)}
$$
for some $0  \neq  \hat{\mathcal{V}}^{\circ}_{\infty}\in \frak g$. 
Let $\hat{\eta}^{\circ}(j)$ and $\hat{\eta}^{\bullet}(j)$ be the Hamiltonian functions
associated to the vector filelds $\hat{\mathcal{V}}^{\circ}(j)$ 
and $\hat{\mathcal{V}}^{\bullet}(j)$ by
$$
\begin{cases}
\; i(\hat{\mathcal{V}}^{\circ}(j))\,\omega (j) \; &=\; \sqrt{-1}\,\bar{\partial}(\hat{\eta}^{\circ}(j)),\\
\; i(\hat{\mathcal{V}}^{\bullet}(j))\,\omega (j) \; &=\; \sqrt{-1}\,\bar{\partial}(\hat{\eta}^{\bullet}(j)),
\end{cases}
$$
where the functions $\hat{\eta}^{\circ}(j)$ and $\hat{\eta}^{\bullet}(j)$
are normalized by the vanishing of the integrals
$\int_M {\hat{\eta}}^{\circ}(j) \omega (j)^n$ and
$\int_M {\hat{\eta}}^{\bullet}(j) \omega (j)^n$. 
Then
$$
\hat{\eta}(j)\; =\; \hat{\eta}^{\circ}(j) + \hat{\eta}^{\bullet}(j).
\leqno{(5.29)}
$$
In view of (5.28), there exists a non-constant real-valued $C^{\infty}$ function 
$\rho$ on $M$
such that $i(\hat{\mathcal{V}}^{\circ}_{\infty})\,\omega^{}_{\infty} \, =\, \sqrt{-1}\,\bar{\partial}\rho $ and that
$$
\hat{\eta}^{\circ}(j) \to \rho \text{ in $C^{\infty}(M)$, as $j \to \infty$.}
$$
Hence by (5.29), it follows from (5.27) 
that
$$
\hat{\eta}^{\bullet}(j) \to - \rho \text{ in $C^{0}(M)$, as $j \to \infty$.}
\leqno{(5.30)}
$$ 
On the other hand, by (5.20), $\int_M |\bar{\partial}\hat{\eta}^{\bullet}(j)|^2_{\omega (j)}\, 
\omega (j)^n \to 0\,$ as $j \to \infty$, and hence 
for each fixed smooth $(0,1)$-form $\theta$ on $M$, we have
\begin{align*}
&\left |\,
(\hat{\eta}^{\bullet}(j), \bar{\partial}(j)_{}^*\theta )^{}_{L^2(M,\, \omega (j)^n )} 
\,\right |
 \; =\;
\left | \int_M (\bar{\partial}\hat{\eta}^{\bullet}(j), \theta )_{\omega(j)}\,
\omega (j)^n \right |\\
&\leq \; \left \{\int_M |\bar{\partial}\hat{\eta}^{\bullet}(j)|^2_{\omega (j)}\, 
\omega (j)^n \right \}^{1/2}\left \{
\int_M |\theta |_{\omega (j)}^2\, \omega (j)^n\right \}^{1/2} \, \to\, 0,
\end{align*}
as $j \to\infty$, where $\bar{\partial}(j)_{}^*$ and  $\bar{\partial}_{\infty}^*$
are the formal adjoints of the operator $\bar{\partial}$ on functions 
for the K\"ahler manifolds $(M, \omega (j))$ and $(M, \omega_{\infty})$, 
respectively.
 Then by letting to $j \to \infty$, we obtain the vanishing 
for the Hermitian $L^2$-inner product of functions $\rho$ and $\bar{\partial}_{}^*\theta$,
$$
(\rho, \bar{\partial}_{}^*\theta )^{}_{L^2(M,\, \omega_{\infty}^n )} \; =\; 0,
$$
for every smooth $(0,1)$-form $\theta$ on $M$, i.e., 
$\bar{\partial} \rho =0$ in a weak sense, and therefore in a strong sense.
Hence we conclude that $\rho$ is constant on $M$ in contradiction to  
$\hat{\mathcal{V}}^{\circ}_{\infty} \neq 0$.

\section{Appendix}

In this appendix, we shall show that the family of K\"ahler manifolds 
$$
(M,\,q\, (\mu_{m,s}^*{\omega_{\operatorname{FS}}})_{|M}),
\qquad -\delta \leq s \leq \delta,\; m=1,2,\dots,
\leqno{(6.1)}
$$
has bounded geometry.  
By Fact 2.6 applied to $m =1$, we identify $\Bbb P^*(E_1)$ with 
$\Bbb A^1 \times \Bbb P^*((E_1)_0) $, and let
$\operatorname{pr}_2:\Bbb P^*(E_1)\to\Bbb P^*((E_1)_0) $ be the 
projection to the second factor. As in Section 4, we have
$$
\mathcal{M}\; \hookrightarrow \; \Bbb P^*(E_1),
\leqno{(6.2)}
$$
where the pullback $\mathcal{H} := 
\operatorname{pr}_2^*\mathcal{O}_{\Bbb P^*((E_1)_0)}(1)$ 
to $\Bbb P^*(E_1)$ of the the hyperplane bundle 
$\mathcal{O}_{\Bbb P^*((E_1)_0)}(1)$
on $\Bbb P^*((E_1)_0)$ is written as
$$
\mathcal{H}_{|\mathcal{M}}\;  = \;  \mathcal{L} \; \, (=\; L^{\ell}).
\leqno{(6.3)}
$$
Recall that the action of $T=\Bbb C^*$ on $\mathcal{M}$ lifts to 
an action of $T$ on $\mathcal{L}$,  and hence $T$ acts on 
$E_1 = \Bbb A^1\times (E_1)_0$ by
$$
T \times (\Bbb A^1 \times (E_1)_0) \to \Bbb A^1 \times (E_1)_0,
\qquad (t, (z, e)) \mapsto (tz, \psi_1 (t)\cdot e).
$$
This induces a $T$-action on $\Bbb P^*(E_1) \, (= \Bbb A^1 \times \Bbb P^*((E_1)_0)) $,
and for (6.2), $\mathcal{M}$ is preserved by the $T$-action.
Note that the $T$-action on $\mathcal{L}$ lifts the
$T$-action on $\mathcal{M}$. By 
$$
T \times \mathcal{M} \to \mathcal{M}, 
\qquad (t, p) \mapsto g_{\mathcal{M}}^{}(t)\cdot p,
$$
we mean the $T$-action on $\mathcal{M}$, and 
the corresponding $T$-action on $\mathcal{L}\otimes \bar{\mathcal{L}}$ 
upstairs will be denoted by
$$
T\, \times\, (\mathcal{L}\otimes \bar{\mathcal{L}}) \to  \mathcal{L}\otimes \bar{\mathcal{L}},
\qquad (t, h) \mapsto g^{}_{|\mathcal{L}|^2} (t)\cdot h.
$$
Since $\operatorname{GL}((E_m)_0)$ acts on $\Bbb P^*((E_m)_0)$
via the projection of
$\operatorname{GL}((E_m)_0)$ onto $\operatorname{PGL}((E_m)_0)$,
by setting $\tilde{\mu}_{m,s} := \psi_m (\exp (s))$, we have 
$$
q\,\mu_{m,s}^*\omega_{\operatorname{FS}} \; =\; 
q\,\tilde{\mu}_{m,s}^*\omega_{\operatorname{FS}}.
\leqno{(6.4)}
$$
In view of
$\delta = C_3 (\log  m) q$, $m \gg 1$, we estimate
$\exp (s)$ in the form
$$
1-\epsilon \leq e^{-C_3 (\log m)/m}\leq \exp(s) 
\leq e^{C_3 (\log m)/m}\leq 1+ \epsilon
\leqno{(6.5)}
$$
for some $0 < \epsilon \ll 1$. 
As in Section 5, by the bases 
 $\{\, \tau_1, \tau_2, \dots, \tau_{N_m}\}$ 
and $\{\, \tau_1', \tau_2', \dots, \tau'_{N_m}\}$ 
for $(E_m)_0$ and $(E_m)_1\, (= V_m) $, respectively,
the spaces $\Bbb P^*((E_m)_0)$ and 
$\Bbb P^*((E_m)_1)\, (= \Bbb P^*(V_m))$ are identified with 
$$
\Bbb P^{N_m-1}(\Bbb C ) =\{\,(z_1: z_2: \dots : z_{N_m})\,\}.
$$
Note that $q \omega_{\operatorname{FS}} = 
(\sqrt{-1}/2\pi ) \partial\bar{\partial}\log \Omega_{\operatorname{FS},m}$,
where $\Omega_{\operatorname{FS},m}$ denotes the positive real smooth section 
$\{(n!/m^n)\Sigma_{\alpha =1}^{N_m}\, |z_{\alpha} |^2 \}^{q}$
of $\mathcal{H}\otimes \bar{\mathcal{H}}$. 
In view of (6.3), identifying $M$ with $\mathcal{M}_1$,
we easily see that $ q\, \tilde{\mu}_{m,s}^*\omega^{}_{\operatorname{FS}}$ 
is 
 $$
(\sqrt{-1}/2\pi ) \, g_{\mathcal{M}}(\exp (s))^*
\partial\bar{\partial}\log \{ g^{}_{|\mathcal{L}|^2}(\exp (s) )\cdot 
\Omega_{\operatorname{FS},m}\},
\leqno{(6.6)}
 $$
when restricted to $M$. By (3.6) and (5.3), we now conclude from (6.4), (6.5) and (6.6)
 that the family of K\"ahler 
manifolds in (6.1) has bounded geometry, as required.



\bigskip\noindent
{\footnotesize
{\sc Department of Mathematics}\newline
{\sc Osaka University} \newline
{\sc Toyonaka, Osaka, 560-0043}\newline
{\sc Japan}}
\end{document}